\documentclass[reqno,11pt]{amsart}
\usepackage{color}
\usepackage{ifpdf}
\ifx\pdftexversion\undefined 
  \usepackage{graphicx}
\else 
  \usepackage[pdftex]{graphicx}
\fi
\ifpdf 
  \DeclareGraphicsExtensions{.pdf,.png,.mps}
  \usepackage{pgf}
\else 
  \usepackage{graphicx}
  \DeclareGraphicsExtensions{.eps,.bmp}
  \usepackage{pstricks}
\fi
\usepackage{floatflt}
\usepackage{wrapfig}
\newcommand{\comment}[1]{}

\makeatletter
 \def\markboth#1#2{%
  \begingroup
   \@temptokena{{#1}{#2}}\xdef\@themark{\the\@temptokena}%
   \mark{\the\@temptokena}%
  \endgroup
  \if@nobreak\ifvmode\nobreak\fi\fi}
  \def\thanks#1{\g@addto@macro\thankses{\thanks{#1}}
  }
\makeatother

\theoremstyle{definition}

\theoremstyle{remark}

\newcommand{\abs}[1]{\left\vert#1\right\vert}
\newcommand{\set}[1]{\left\{#1\right\}}

\def\deg{\mathop {\rm deg}}
\def\indeg{\mathop {\rm indeg}}
\def\lead{\mathop {\rm lead}}

\begin{document}

\title[Reverse Pr\"ufer Algorithm]%
{A Generalized Enumeration of Labeled Trees and Reverse Pr\"ufer Algorithm}%
\author[S.~Seo]{Seunghyun Seo}%
\address{Department of Mathematics\\
Seoul National University\\
Seoul, 151-742\\
Republic of Korea}%
\email{shseo@snu.ac.kr}%

\author[H.~Shin]{Heesung Shin}%
\address{Department of Mathematics\\
Korea Advanced Institute of Science and Technology\\
Daejeon, 305-701\\
Republic of Korea}%
\email{H.Shin@kaist.ac.kr}%

\today


\begin{abstract}
A {\em leader} of a tree $T$ on $[n]$ is a vertex which has no
smaller descendants in $T$. Gessel and Seo showed
$$\sum_{T \in \mathcal{T}_n}u^\text{(\# of leaders in $T$)}\,c^\text{(degree of $1$ in $T$)}=u\,P_{n-1}(1,u,cu),$$
which is a generalization of Cayley formula, where $\mathcal{T}_n$
is the set of trees on $[n]$ and
$$P_n(a,b,c)=c\prod_{i=1}^{n-1}(ia+(n-i)b+c).$$ Using
a variation of Pr\"ufer code which is called a {\em RP-code}, we
give a simple bijective proof of Gessel and Seo's formula.
\end{abstract}
\maketitle
\section{Introduction}
A {\em tree} on $V$ is an acyclic connected graph with vertex set
$V$. In 1889, Cayley~\cite{C89} showed that
$\abs{\mathcal{T}_n}=n^{n-2}$ ($n\ge 1$), called {\em Cayley
formula}, where $\mathcal{T}_n$ is the set of trees on
$[n]=\set{1,2,\ldots,n}$. Later, in 1918, Pr\"ufer~\cite{P18} made
the {\em Pr\"ufer code} which is a bijection between $\mathcal{T}_n$
and $[n]^{n-2}$. Assume that edges are directed toward the vertex
$1$ and $\indeg_T(i)$ is the indegree of $i$ in $T$. By Pr\"ufer
code, we have
$$\sum_{T \in \mathcal{T}_n} \prod_{i \in [n]} x_i^{\indeg_T(i)}=
x_1(x_1+\cdots + x_n)^{n-2},$$ which is a generalization of Cayley
formula.

A tree is called a {\em rooted tree} if one vertex has been
designated the root. A vertex $v$ in a rooted tree is a {\em
descendant} of $u$ if $u$ lies on the unique path from the root to
$v$. By convention, we consider that (unrooted) trees are rooted at
the smallest vertex. A vertex $v$ of a rooted tree is called a {\em
leader} if $v$ is minimal among its descendants. Note that `leader'
is the new terminology of `proper vertex' which was introduced by
Seo \cite{S04}.

Recently, Gessel and Seo~\cite{GS05} showed that
\begin{equation} \label{main}
\sum_{T \in \mathcal{T}_n} u^{\lead(T)} c^{\deg_T(1)} = u\, P_{n-1}
(1,u,cu),
\end{equation}
where $\lead(T)$ is the number of leaders in $T$ and the homogeneous
polynomial $P_n(a,b,c)$ is defined by
$$P_n(a,b,c)=c\prod_{i=1}^{n-1}(ia+(n-i)b+c).$$
To prove \eqref{main}, they used generating functions methods.

In this paper, we prove the equation \eqref{main} by giving an
algorithm which produces a code with length $n-1$ from a tree with
$n$ vertices.

\section{Reverse Pr\"ufer Algorithm}
The {\em reverse Pr\"ufer code (RP-code)} $\varphi(T)=(
\sigma_{1},\ldots,\sigma_{n-1})$ of a rooted tree $T$ on $[n]$ is
generated by successively selecting the unselected vertex of $T$
having the smallest descendants including itself. If several
vertices have the same smallest descendant, we choose the vertex
which is the closest to the root. Because the root is selected above
all, we assume that the root was already selected. To obtain the
code from $T$, we select such a vertex in each step, recording its
parent $\sigma_i$, from the tree, until all the vertices are
selected. We call this process a {\em reverse Pr\"ufer algorithm
(RP-algorithm)}.

\begin{figure}[t]
\input{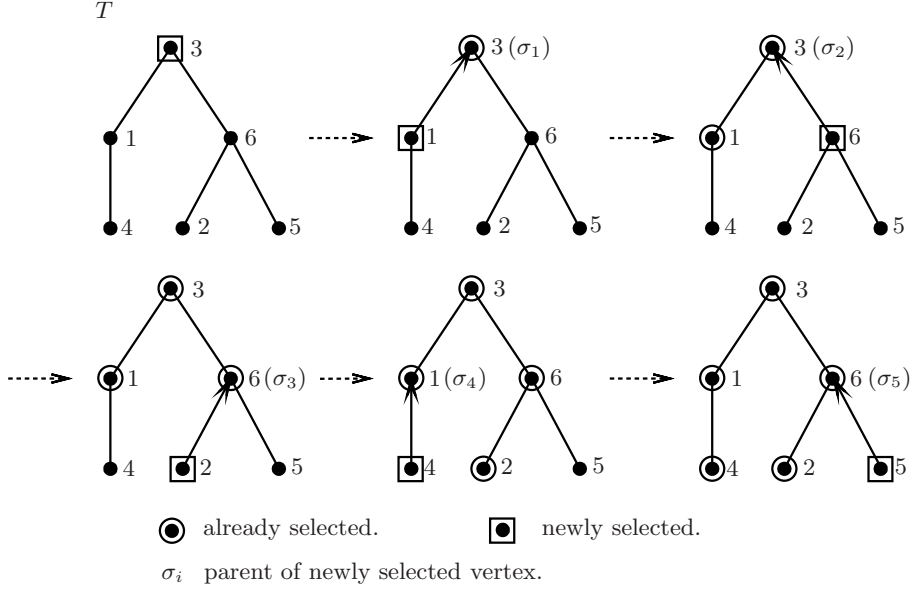}
\caption{Rooted Tree $T$ to RP-code $\varphi(T)=(3,3,6,1,6)$}
\label{tree2code}
\end{figure}

The inverse of $\varphi$ is described as follows: Let
$\sigma=(\sigma_{1},\dots,\sigma_{n-1})$ be a sequence of positive
integers with $\sigma_i\in [n]$ for all $i$. We can find the tree
$T$ whose code is $\sigma$ as building up labeled trees $T_i$ with
$i+1$ vertices, except one leaf is unlabeled, by reading the code
$\sigma$ forward. Before reading the code, we consider the rooted
tree $T_0$ with only one vertex. This root of $T_0$ is unlabeled.
Assume that $T_{i-1}$ is the labeled tree which corresponds to
initial $i-1$ code $(\sigma_1,\ldots,\sigma_{i-1})$ for
$i=1,\ldots,n-1$. We make $T_i$ as follows: We label $\sigma_i$ to
the unlabeled leaf of $T_{i-1}$. But if $\sigma_i$ is belong to
labels of $T_{i-1}$, use the minimum of unused labels instead of
$\sigma_i$ as new label number. And then add an unlabeled vertex and
an edge between $\sigma_i$ and the just added vertex. After reading
the code $\sigma$, we obtain $T_{n-1}$ with $n$ vertices. The
unlabeled vertex of $T_{n-1}$ is labeled by the unused label among
$[n]$. Then we get the tree $T$. Note that the map $\varphi^{-1}$
was already mentioned in \cite[pp. 1--2]{CKSS04}.

\begin{figure}[t]
\input{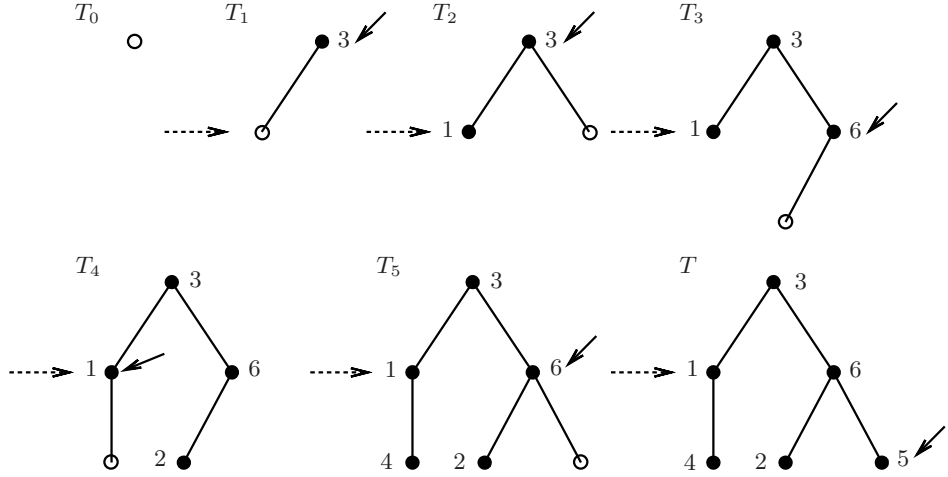}
\caption{RP-code $\sigma=(3,3,6,1,6)$ to Tree
$\varphi^{-1}(\sigma)=T$} \label{code2tree}
\end{figure}

The first coordinate of the RP-code of a rooted tree $T$ is always
the label of the root of $T$. In particular, the RP-code of a tree
on $[n]$ begins with~1. Cayley formula is reconfirmed by the number
of RP-codes.

\begin{figure}[t]
\input{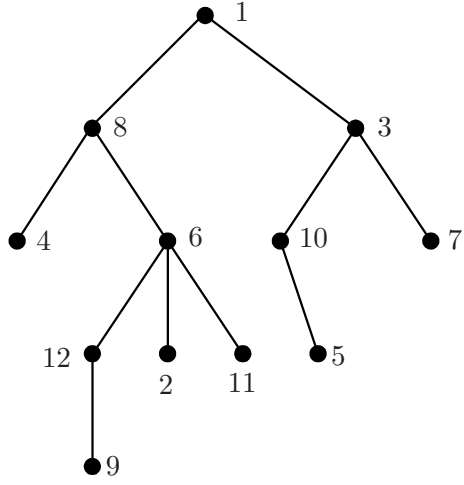}
\caption{Example of the tree with root 1} \label{example}
\end{figure}
Figure~\ref{example} shows the tree corresponding to the RP-code
$(1,8,6,1,8,3,10,\linebreak[2]3,6,12,6)$.

\section{Statistics of Leader in Trees}
Now we trace leaders in $T$ during the RP-alogrithm. Let
$\sigma=(\sigma_1,\ldots,\sigma_n)$ be a RP-code. For each
$i=2,\ldots,n$, let $T_{i-1}$ is the tree obtained from subcode
$\sigma_1,\ldots,\sigma_{i-1}$. Let $l$ be a minimal element in
$[n]$ which does not appear in $T_{i-1}$. To construct $T_i$ from
$T_{i-1}$ and $\sigma_i$, we should consider the following two
cases.
\begin{enumerate}
\item \label{case1} Suppose that $\sigma_i$ appears in $T_{i-1}$.
Then the unlabeled vertex $v$ in $T_{i-1}$ is labeled by $l$ in
$T_i$. Since the new label $l$ is minimal among unused labels in
$T_{i-1}$, the vertex $v$ is a leader in $T$.


\item Suppose that $\sigma_i$ does not appear in $T_{i-1}$.
Then the unlabeled vertex $v$ in $T_{i-1}$ is labeled by $\sigma_i$
in $T_i$.
\begin{enumerate}
\item \label{case2a} If $\sigma_i=l$, then the vertex $v$ is leader in $T$ like case~(\ref{case1}).
\item \label{case2b} If $\sigma_i\ne l$, then the vertex $v$ has a descendent labeled by
$l$. Thus, the vertex $v$ is not leader in $T$.
\end{enumerate}
\end{enumerate}

So there are exactly $i$ choices of $\sigma_i$, case~(\ref{case1})
and case~(\ref{case2a}), such that the newly labeled vertex $v$ is a
leader in $T$. Because the number of $r$'s (=$\sigma_1$) in a
RP-code equals to the degree of the root $r$ in $T$, $\deg_T(1)$ is
the number of $1$ in the RP-code of a tree $T$.

Thus we have the following formula:

\begin{alignat*}{2}
\sum_{T \in \mathcal{T}_n} u^{\lead(T)} c^{\deg_T(1)}&&~=~& cu \tag*{by $\sigma_1(=1)$}\\
&&&\times ((n-2)+u+cu) \tag*{by $\sigma_2$}\\
&&&\times ((n-3)+ 2u +cu) \tag*{by $\sigma_3$}\\
&&&\qquad \vdots\\
&&&\times (1 + (n-2)u + cu) \tag*{by $\sigma_{n-1}$}\\
&&&\times u \tag*{by filling the last label}\\
&&\,=\,& c u^2 \prod_{i=2}^{n-1}((n-i)+(i-1)u+cu).\\
&&\,=\,& u P_{n-1} (1,u,cu).
\end{alignat*}
This completes the bijective proof of equation~\eqref{main}.
\section{Remarks}

\begin{enumerate}
\item If $(a_1,\ldots,a_{n-2},1)$ is a Pr\"ufer code of $T$ and
$\varphi(T)=(1,\sigma_2,\ldots,\sigma_{n-1})$ is a RP-code of $T$,
then $a_i=\sigma_{n-i}$ for each $i$. This justifies the terminology
{\em `reverse' Pr\"ufer code.}
\item With a slight variation of the RP-algorithm, we also find a
combinatorial proof of the following formulas for $k$-ary trees and
ordered trees.
\begin{eqnarray*}
\sum_U u^{\lead(U)} &=&
P_n(k,(k-1)u,u)\\
\sum_V u^{\lead(V)} &=&
P_n(1,2u,u)\\
\end{eqnarray*}
where $U$ runs all $k$-ary trees and $V$ runs all ordered trees on
$[n]$.

\end{enumerate}

\bibliographystyle{amsplain}

\begin{thebibliography}{AAA}
\bibitem{C89}
A.~Cayley, {\em A theorem on trees}, Quart. J. Math. {\bf 23},
376--378, 1889.

\bibitem{P18}
H.~Pr\"ufer, {\em Neuer Beweis eines Satzes \"uber Permutationen},
Archiv der Math. (3) {\bf 27}, 142--144, 1918.


\bibitem{S04}
S.~Seo, {\em A combinatorial proof of Postnikov's identity and a
generalized enumeration of labeled trees}, Electron. J. Combin. {\bf
11} no. 2, \#N3, 2004.

\bibitem{CKSS04}
M.~Cho, D.~Kim, S.~Seo and H.~Shin, {\em Colored Prufer Codes for
k-Edge Colored Trees}, Electron. J. Combin. {\bf 11} no. 1, \#N10,
2004.

\bibitem{GS05}
Ira~M.~Gessel and S.~Seo, {\em A refinement of Cayley's formula for
trees}, arXiv: math.CO/0507497, 2005.
\end{thebibliography}

\end{document}